\documentclass[twoside,reqno,A4]{amsart}

\usepackage{amsmath}
\newcommand{\qdn}{\hspace*{-1.5mm}}
\newcommand{\qqdn}{\hspace*{-2.5mm}}

\newcommand{\+}{&\qqdn}%

\newcommand{\ffnk}[4]{\left[\qdn\ba{#1}#3\\#4\ea{\!\Big|\:#2}\right]}

\newcommand{\nnm}{\nonumber}
\newcommand{\be}{\begin{equation}}
\newcommand{\ee}{\end{equation}}
\newcommand{\ba}{\begin{array}}
\newcommand{\ea}{\end{array}}
\newcommand{\bmn}{\begin{eqnarray}}
\newcommand{\emn}{\end{eqnarray}}
\newcommand{\bnm}{\begin{eqnarray*}}
\newcommand{\enm}{\end{eqnarray*}}
\newcommand{\bln}{\begin{subequations}}
\newcommand{\eln}{\end{subequations}}

\newtheorem{thm}{Theorem}[section]

\newtheorem{prop}[thm]{Proposition}

\newcommand{\bbtm}[4]{\bibitem{kn:#1}{#2,}~{#3,}~{#4.}}
\newcommand{\cito}[1]{\cite{kn:#1}}

\begin{document}

\title{\large Infinite Summation Formulas Involving Riemann-Zeta function}
\author{Xiaoxia Wang, Xueying Yuan$^{\dag}$}
\dedicatory{Department of Mathematics,
            Shanghai University,\\
            Shanghai 200444, P.\:R.\:China}

\thanks{$^{\dag}$ Corresponding author. \\
E-mail addresses: yuanxueying@shu.edu.cn (X. Yuan), xiaoxiawang@shu.edu.cn (X. Wang).\\
This work is supported by National Natural Science Foundations of China (11661032).}

\maketitle\thispagestyle{empty}
\markboth{X. Wang-X. Yuan}%
       {Infinite Summation Formulas related Riemann-Zeta function}
\begin{center}\parbox{110mm}{

By some hypergeometric summation theorems, the authors establish a series of new infinite summation formulas involving generalized harmonic numbers related to Riemann-Zeta function, with three different patterns.\\

}
\end{center}
\vspace*{10mm}
\section{Introduction}
Following Slater\cito{Slater}, the generalized hypergeometric series is defined by
\bnm \label{def-de}
\+\+ {_{p+1}F_q}
\ffnk{ccc}{z}{\alpha_0,\alpha_1,\+\cdots,\+\alpha_p}
              {\beta_1,\+\cdots,\+\beta_q}
:=\sum_{n=0}^\infty
\frac{(\alpha_0)(\alpha_1)_n\cdots(\alpha_p)_n}
     {(\beta_1)_n\cdots(\beta_q)_n}
 \frac{z^n}{n!},
 \enm
where the shifted factorial is defined by
\bnm
(x)_0:\equiv 1 \quad \text{and} \quad (x)_n:=\frac{\Gamma(x+n)}{\Gamma(x)}=x(x+1)\cdots(x+n-1),\quad n=1,2,\cdots.
\enm
For more details of hypergeometric series, we refer the readers to \cito{Slater}.
For the $\Gamma$-function, there hold the Weierstrass product expression \cito{Andrews}
\bnm
\Gamma(z)=z^{-1}\prod\limits_{n=1}^\infty{(1+1/n)^z/(1+z/n)},
\enm
and the logarithmic differentiation with the Euler constant being given by
\bnm
\frac{\Gamma'(z)}{\Gamma(z)}=-\gamma+\sum\limits_{n=0}^\infty\frac{z-1}{(n+1)(n+z)}, \quad\gamma=\lim\limits_{n\rightarrow \infty}\left\{\sum\limits_{k=1}^n\frac{1}{k}-\ln n\right\}.
\enm
Also the following expansions of $\Gamma$-function are found in \cito{Zheng},
\bmn
\+\+\Gamma(1-z)=\exp\left\{\sum\limits_{k=1}^\infty\frac{\sigma_k}{k}z^k\right\};\label{1.1}\\[2mm]
\+\+\Gamma(\frac{1}{2}-z)=\sqrt{\pi}\exp\left\{\sum\limits_{k=1}^\infty\frac{\tau_k}{k}z^k\right\},\label{1.2}
\emn
where the Riemann-Zeta sequences $\{\sigma_k,\tau_k\}$ are defined by
\bnm
\+\+\sigma_1=\gamma,\quad\quad\sigma_m=\zeta(m)=\sum\limits_{k=1}^\infty\frac{1}{k^m},\quad m=2,3,\cdots;\\[2mm]
\+\+\tau_1=\gamma+2\ln2,\quad\quad\tau_m=(2^m-1)\zeta(m),\quad m=2,3,\cdots.
\enm
Throughout the paper, the Euler summation formulas $\zeta(2)$ and $\zeta(4)$ are expressed by $\pi$ as follows \cito{Chu}:
\bnm
\zeta(2)=\sum\limits_{n=1}^\infty\frac{1}{n^2}=\frac{\pi^2}{6},\qquad\zeta(4)=\sum\limits_{n=1}^\infty\frac{1}{n^4}=\frac{\pi^4}{90}.
\enm
The generalized harmonic numbers which are defined by
\bnm
\+\+H_n=\sum\limits_{k=1}^n\frac{1}{k},\quad\quad H_n^{(r)}=\sum\limits_{k=1}^n\frac{1}{k^r},\quad r=2,3,\cdots,\\[2mm]
\+\+O_n=\sum\limits_{k=1}^n\frac{1}{2k-1},\quad\quad O_n^{(r)}=\sum\limits_{k=1}^n\frac{1}{(2k-1)^r},\quad r=2,3,\cdots,
\enm
can be used to express the following finite products \cito{Macdonald} through the symmetric functions:
\bmn
\+\+\prod\limits_{k=1}^n(1+\frac{x}{k})=1+xH_n+\frac{x^2}{2}(H_n^2-H_n^{(2)})+\frac{x^3}{6}(H_n^3-3H_n H_n^{(2)}+2H_n^{(3)})+\cdots;\label{2}\\
\+\+\prod\limits_{k=1}^n(1-\frac{x}{k})^{-1}=1+xH_n+\frac{x^2}{2}(H_n^2+H_n^{(2)})+\frac{x^3}{6}(H_n^3+3H_n H_n^{(2)}+2H_n^{(3)})+\cdots;\label{3}\\
\+\+\prod\limits_{k=1}^n(1+\frac{y}{2k-1})=1+yO_n+\frac{y^2}{2}(O_n^2-O_n^{(2)})+\frac{y^3}{6}(O_n^3-3O_n O_n^{(2)}+2O_n^{(3)})+\cdots;\label{4}\\
\+\+\prod\limits_{k=1}^n(1-\frac{y}{2k-1})^{-1}=1+yO_n+\frac{y^2}{2}(O_n^2+O_n^{(2)})+\frac{y^3}{6}(O_n^3+3O_n O_n^{(2)}+2O_n^{(3)})+\cdots.\label{5}
\emn

The identities involving generalized harmonic numbers related to Riemann-Zeta function can be established by many different methods. For example, De Doelder \cito{Doelder} applied the irregular integral and Borwein \cito{Borwein} used the Parseval identity on Fourier series to get some of the identities. Also, the hypergeometric method \cito{Shen} is used to obtain infinite summation formulas involving generalized harmonic numbers related to Riemann-Zeta function. More interesting infinite summation formulas related to Riemann-Zeta function have been found, such as, in \cite{{kn:Chen, kn:Chu, kn:ZW, kn:WC}}.\\

Applying the hypergeometric method, we will establish some infinite summation formulas involving generalized harmonic numbers related to Riemann-Zeta function from the famous summation theorems due to Gauss, Watson and Bailey. Here, we will present the following three different patterns which have not been found before
\bnm
\+\+\sum_{k=1}^\infty\frac{\binom{2k}{k}}{k^i2^{2k}}P_k, \qquad i=1,2,3;\\
\+\+\sum_{k=1}^\infty\frac{P_k}{k^i2^k},\quad i=1,2;\\
\+\+\sum_{k=1}^\infty\frac{3^k}{k^2\binom{2k}{k}}P_k,
\enm
where $P_k$ is a polynomial in $H_k^{(r)}$ or $O_k^{(r)}$$\:{(k,r\in Z^+)}$.

\section{Summation Formulas Involving  Riemann-Zeta function from Gauss and Watson summation theorems}
In this section, we will derive some infinite summation formulas related to Riemann-Zeta function from the Gauss summation theorem \eqref{GG1} and Watson
summation theorem \eqref{2-11} with the pattern as follows:
\bnm
\sum_{k=1}^\infty\frac{\binom{2k}{k}}{k^i2^{2k}}P_k, \qquad i=1,2,3,
\enm
where $P_k$ is a polynomial in $H_k^{(r)}$ or $O_k^{(r)}$$\:{(k,r\in Z^+)}$.\\

\begin{thm}[Gauss, \cito{Slater} ] For complex parameters $a,b,c$ with $\mathcal{R}(c-a-b)>0$, the following summation formula is true.
\bmn\label{GG1}
{_2F_1}\ffnk{cc}{1}{a,\+b}{\+c}
=\frac{\Gamma(c)\Gamma(c-a-b)}{\Gamma(c-a)\Gamma(c-b)}
\emn
\end{thm}
Performing the replacements $a\rightarrow a+\frac{1}{2}$ and $c\rightarrow c+1$ in \eqref{GG1}, we attain the following expression:
\bmn
{_2F_1}\ffnk{c}{1}{\frac{1}{2}+a,b}{1+c}
=\frac{\Gamma(1+c)\Gamma(\frac{1}{2}+c-a-b)}{\Gamma(\frac{1}{2}+c-a)\Gamma(1+c-b)}\label{G1}
\emn
Recalling the definition of hypergeometric series and applying \eqref{1.1} and \eqref{1.2}, we can write \eqref{G1} as
\bmn{\label{c11}}
\+\+\quad1+b\sum\limits_{k=1}^\infty\frac{\binom{2k}{k}\prod\limits_{i=1}^{k}(1+\frac{2a}{2i-1})\prod\limits_{m=1}^{k-1}(1+\frac{b}{m})}{k2^{2k}\prod\limits_{j=1}^{k}(1+\frac{c}{j})}\\
\+\+\quad=\exp\Big\{\sum\limits_{k=1}^\infty\frac{(-1)^k\tau_k}{k}\big[(c-a-b)^k-(c-a)^k\big]
+\sum\limits_{k=1}^\infty\frac{(-1)^k\sigma_k}{k}\big[c^k-(c-b)^k\big]
\Big\}\\
\+\+\quad=\exp\Big\{\tau_1b
+\tau_2(ab-bc+\frac{1}{2}b^2)
+\tau_3(a^2b+ab^2-2abc+\frac{1}{3}b^3-b^2c+bc^2)\nnm
\emn
\bmn
\+\+\quad\quad+\tau_4(a^3b+\frac{3}{2}a^2b^2-3a^2bc+ab^3-3ab^2c+3abc^2+\frac{1}{4}b^4-b^3c+\frac{3}{2}b^2c^2-bc^3)
\nnm\\
\+\+\quad\quad+\tau_5(a^4b+2a^3b^2-4a^3bc+2a^2b^3-6a^2b^2c+ 6a^2bc^2+ab^4-4ab^3c+6ab^2c^2\nnm\\
\+\+\quad\quad-4abc^3+\frac{1}{5}b^5-b^4c+2b^3c^2-2b^2c^3+bc^4)
+\cdots
\Big\}\times\exp\Big\{-\sigma_1b
+\sigma_2(-\frac{1}{2}b^2
\nnm\\
\+\+\quad\quad
+bc)+\sigma_3(-\frac{1}{3}b^3+b^2c-bc^2)+\sigma_4(-\frac{1}{4}b^4+b^3c-\frac{3}{2}b^2c^2+bc^3)
+\sigma_5(-\frac{1}{5}b^5\nnm\\
\+\+\quad\quad+b^4c-2b^3c^2+2b^2c^3-bc^4)+\cdots
\Big\}.\label{2-10}
\emn
Expanding the infinite summation over products in \eqref{c11} to obtain another alternative multivariate series expansion with generalized harmonic number by employing relations \eqref{2}--\eqref{4}. Applying the relations (5) and (6), we can obtain a compact writing of a multivariate series by the power series expansions in \eqref{2-10} with Riemann-Zeta function. Therefore, comparing the coefficients of the above two multivariate series expansions term-by-term, we get a number of infinite summation formulas involving generalized harmonic numbers related to Riemann-Zeta function. Here, we use a self-explanatory notation $[a^ib^jc^k]$ to show the process of extracting the coefficients of monomial $a^ib^jc^k$ from multivariate power series expansions.
\begin{prop}[Infinite summation formulas related to $\zeta(2)$]\label{exam2.3}
\bmn
\+\+[ab]\quad\quad\sum\limits_{k=1}^\infty\frac{\binom{2k}{k}O_{k}}{k2^{2k}}=\frac{\pi^2}{4};\label{g1.1.1}\\
\+\+[bc]\quad\quad\sum\limits_{k=1}^\infty\frac{\binom{2k}{k}H_{k}}{k2^{2k}}=\frac{\pi^2}{3};\label{g1.1.2}\\
\+\+[b^2]\quad\quad\sum\limits_{k=1}^\infty\frac{\binom{2k}{k}H_{k-1}}{k2^{2k}}=\frac{\pi^2}{6}+2{\ln^22}.\label{g1.1.5}
\emn
\end{prop}
\begin{prop}[Infinite summation formulas related to $\zeta(3)$]
\bmn
\+\+[a^2b]\quad\quad\sum\limits_{k=1}^\infty\frac{\binom{2k}{k}(O^2_{k}-O_k^{(2)})}{k2^{2k}}=\frac{7\zeta(3)}{2};\label{g1.8}\\
\+\+[bc^2]\quad\quad\sum\limits_{k=1}^\infty\frac{\binom{2k}{k}(H^2_{k}+H_k^{(2)})}{k2^{2k}}=12\zeta(3);\label{g1.11}\\
\+\+[abc]\quad\quad\sum\limits_{k=1}^\infty\frac{\binom{2k}{k}{H_k}O_k}{k2^{2k}}=7\zeta(3)\label{g1.1}.
\emn
\end{prop}

\begin{prop}[Infinite summation formulas related to $\zeta(2)$ and $\zeta(3)$]
\bmn
\+\+[ab^2]\quad\quad\sum\limits_{k=1}^\infty\frac{\binom{2k}{k}H_{k-1}{O_k}}{k2^{2k}}=\frac{7}{2}\zeta(3)+\frac{\pi^2}{2}(\ln2);\label{g1.5}\\
\+\+[b^2c]\quad\quad\:\sum\limits_{k=1}^\infty\frac{\binom{2k}{k}{H_{k-1}}H_k}{k2^{2k}}=6\zeta(3)+\frac{2\pi^2}{3}(\ln2);\label{g1.6}\\
\+\+[b^3]\quad\quad\:\:\sum\limits_{k=1}^\infty\frac{\binom{2k}{k}(H^2_{k-1}-H_{k-1}^{(2)})}{k2^{2k}}=4\zeta(3)+\frac{2\pi^2}{3}(\ln2)+\frac{8}{3}\ln^32.\label{g1.7}
\emn
\end{prop}
Now we can obtain the following infinite summation formulas $\sum\limits_{k=1}^\infty\frac{\binom{2k}{k}O_k}{k^i2^{2k}}$ with $i=1, 2$.
\bnm
\+\+\eqref{g1.1.1}\quad\quad\quad\quad\sum\limits_{k=1}^\infty\frac{\binom{2k}{k}O_{k}}{k2^{2k}}=\frac{\pi^2}{4};\\
\+\+\eqref{g1.1}-\eqref{g1.5}\quad\sum\limits_{k=1}^\infty\frac{\binom{2k}{k}O_{k}}{k^22^{2k}}=\frac{7\pi^2}{12}-\frac{\pi^2}{2}(\ln2).
\enm

\begin{prop}[Infinite summation formulas related to $\zeta(4)$]
\bmn
\+\+[a^2bc]\quad\:\:\sum\limits_{k=1}^\infty\frac{\binom{2k}{k}H_k(O_k^2-O_k^{(2)})}{k2^{2k}}=\frac{\pi^4}{4};\label{g1.9}\\
\+\+[abc^2]\quad\sum\limits_{k=1}^\infty\frac{\binom{2k}{k}O_k(H_k^2+H_k^{(2)})}{k2^{2k}}=\frac{\pi^4}{2};\\
\+\+[a^3b]\quad\sum\limits_{k=1}^\infty\frac{\binom{2k}{k}(O_k^3-3{O_k}O_k^{(2)}+2O_k^{(3)})}{k2^{2k}}=\frac{\pi^4}{8};\label{g1.14}\\
\+\+[bc^3]\quad\:\:\:\sum\limits_{k=1}^\infty\frac{\binom{2k}{k}(H_k^3+3{H_k}H_k^{(2)}+2H_k^{(3)})}{k2^{2k}}=\frac{14\pi^4}{15}.\label{g1.20}
\emn
\end{prop}

\begin{prop}[Infinite summation formulas related to $\zeta(2)$, $\zeta(3)$ and $\zeta(4)$]
\bmn
\+\+[a^2b^2]\quad\sum\limits_{k=1}^\infty\frac{\binom{2k}{k}H_{k-1}(O_k^{2}-O_k^{(2)})}{k2^{2k}}=\frac{3\pi^4}{16}+7(\ln2)\zeta(3);\label{g1.10}\\
\+\+[b^2c^2]\quad\sum\limits_{k=1}^\infty\frac{\binom{2k}{k}H_{k-1}(H_k^{2}+H_k^{(2)})}{k2^{2k}}=\frac{26\pi^4}{45}+24(\ln2)\zeta(3);\label{g1.13}\\
\+\+[ab^2c]\quad\sum\limits_{k=1}^\infty\frac{\binom{2k}{k}{H_{k-1}}{H_k}O_k}{k2^{2k}}=\frac{\pi^4}{3}+14(\ln2)\zeta(3);\label{g1.15}\\
\+\+[ab^3]\quad\sum\limits_{k=1}^\infty\frac{\binom{2k}{k}O_{k}(H_{k-1}^{2}-H_{k-1}^{(2)})}{k2^{2k}}=\frac{\pi^4}{4}+14(\ln2)\zeta(3)+\pi^2(\ln^2{2});\label{g1.17}\\
\+\+[b^3c]\quad\sum\limits_{k=1}^\infty\frac{\binom{2k}{k}H_{k}(H_{k-1}^{2}-H_{k-1}^{(2)})}{k2^{2k}}=\frac{19\pi^4}{45}+24(\ln2)\zeta(3)+\frac{4\pi^2}{3}(\ln^2{2}).\label{g1.19}
\emn
\end{prop}

\begin{prop}[Infinite summation formulas related to $\zeta(5)$]
\bmn
\+\+[a^3bc]\quad\sum\limits_{k=1}^\infty\frac{\binom{2k}{k}H_{k}(O_k^{3}-3{O_k}O_k^{(2)}+2O_k^{(3)})}{k2^{2k}}=93\zeta(5);\label{g1.1.3}\\
\+\+[abc^3]\quad\sum\limits_{k=1}^\infty\frac{\binom{2k}{k}O_{k}(H_k^{3}+3{H_k}H_k^{(2)}+2H_k^{(3)})}{k2^{2k}}=372\zeta(5);\\
\+\+[a^2bc^2]\quad\sum\limits_{k=1}^\infty\frac{\binom{2k}{k}(O_k^{2}-O_k^{(2)})(H_k^{2}+H_k^{(2)})}{k2^{2k}}=186\zeta(5).\label{g1.1.6}
\emn
\end{prop}

\begin{prop}[Infinite summation formulas to $\zeta(2)$, $\zeta(3)$, $\zeta(4)$ and $\zeta(5)$]
\bmn
\+\+[a^2b^2c]\quad\sum\limits_{k=1}^\infty\frac{\binom{2k}{k}{H_k}H_{k-1}(O_k^2-O_k^{(2)})}{k2^{2k}}=93\zeta(5)+\frac{\pi^4}{2}(\ln2)+\frac{14\pi^2}{3}\zeta(3);\label{g10}\\
\+\+[ab^2c^2]\quad\sum\limits_{k=1}^\infty\frac{\binom{2k}{k}{O_k}H_{k-1}(H_k^2+H_k^{(2)})}{k2^{2k}}=186\zeta(5)+\pi^4(\ln2)+\frac{23\pi^2}{3}\zeta(3);\\
\+\+[a^3b^2]\quad\sum\limits_{k=1}^\infty\frac{\binom{2k}{k}H_{k-1}(O_k^{3}-3{O_k}O_k^{(2)}+2O_k^{(3)})}{k2^{2k}}
=\frac{93}{2}\zeta(5)+\frac{\pi^4}{4}(\ln2)+\frac{21\pi^2}{8}\zeta(3);\label{g1.1.4}\\
\+\+[b^2c^3]\quad\sum\limits_{k=1}^\infty\frac{\binom{2k}{k}H_{k-1}(H_k^{3}+3{H_k}H_k^{(2)}+2H_k^{(3)})}{k2^{2k}}
=360\zeta(5)+\frac{28\pi^4}{15}(\ln2)+12\pi^2\zeta(3);\\
\+\+[a^2b^3]\quad\sum\limits_{k=1}^\infty\frac{\binom{2k}{k}(O_k^2-O_k^{(2)})(H_{k-1}^2-H_{k-1}^{(2)})}{k2^{2k}}\nnm\\
\+\+\quad\quad\quad\quad\quad\quad=62\zeta(5)+\frac{3\pi^4}{4}(\ln2)+\frac{14\pi^2}{3}\zeta(3)+14(\ln^22)\zeta(3);\label{c9}\\
\+\+[b^3c^2]\quad\sum\limits_{k=1}^\infty\frac{\binom{2k}{k}(H_k^2+H_k^{(2)})(H_{k-1}^2-H_{k-1}^{(2)})}{k2^{2k}}\nnm\\
\+\+\quad\quad\quad\quad\quad\quad=240\zeta(5)+\frac{104\pi^4}{45}(\ln2)+12\pi^2\zeta(3)+48(\ln^22)\zeta(3);\\
\+\+[ab^3c]\quad\sum\limits_{k=1}^\infty\frac{\binom{2k}{k}{H_k}O_k(H_{k-1}^2-H_{k-1}^{(2)})}{k2^{2k}}\nnm\\
\+\+\quad\quad\quad\quad\quad\quad=124\zeta(5)+\frac{4\pi^4}{3}(\ln2)+\frac{23\pi^2}{3}\zeta(3)+28(\ln^22)\zeta(3);\\
\+\+[ab^4]\quad\sum\limits_{k=1}^\infty\frac{\binom{2k}{k}O_k(H_{k-1}^3-3{H_{k-1}}H_{k-1}^{(2)}+2H_{k-1}^{(3)})}{k2^{2k}}=93\zeta(5)+\frac{3\pi^4}{2}(\ln2)+\frac{13\pi^2}{2}\zeta(3)\nnm\\
\+\+\quad\quad\quad\quad\quad\quad+42(\ln^22)\zeta(3)+2\pi^2(\ln^32).
\emn
\end{prop}
Therefore, we get the following infinite summation formula.
\bnm
\+\+\eqref{g1.1.3}-\eqref{g1.1.4}\quad\sum\limits_{k=1}^\infty\frac{\binom{2k}{k}(O_{k}^3-3{O_k}O_k^{(2)}+2O_k^{(3)})}{k^22^{2k}}=\frac{93\zeta(5)}{2}-\frac{\pi^4}{4}(\ln2)-\frac{21\pi^2}{8}\zeta(3).
\enm
Applying the formula \eqref{g10}, we obtain the following infinite summation formulas of the pattern $\sum\limits_{k=1}^\infty\frac{\binom{2k}{k}(O_k^2-O_k^{(2)})}{k^i2^{2k}}$ with $i=1, 2, 3$.
\bnm
\+\+\eqref{g1.8}\quad\quad\quad\quad\sum\limits_{k=1}^\infty\frac{\binom{2k}{k}(O^2_{k}-O_k^{(2)})}{k2^{2k}}=\frac{7\zeta(3)}{2};\\
\+\+\eqref{g1.9}-\eqref{g1.10}\quad\sum\limits_{k=1}^\infty\frac{\binom{2k}{k}(O_{k}^2-O_k^{(2)})}{k^22^{2k}}=\frac{\pi^4}{8}-7(\ln2)\zeta(3)-\frac{3\pi^2}{8};\\
\+\+\eqref{g1.1.6}+\eqref{c9}\quad\sum\limits_{k=1}^\infty\frac{\binom{2k}{k}(O_{k}^2-O_k^{(2)})}{k^32^{2k}}=31\zeta(5)-\frac{\pi^4}{8}(\ln2)-\frac{7\pi^2}{3}\zeta(3)+7(\ln^22)\zeta(3).
\enm
Next, we will deduce some infinite summation formulas related to Riemann-Zeta function from the Watson summation theorem.
\begin{thm}[Watson \cito{Slater}] For complex parameters $a,b,c$ with $\mathcal{R}(1-a-b+2c)>0$, the following summation formula is true.
\bmn
{_3F_2}\ffnk{cc}{1}{a,b,c}{\frac{1}{2}+\frac{1}{2}a+\frac{1}{2}b,2c}
=\frac{\Gamma(\frac{1}{2})\Gamma(\frac{1}{2}+c)\Gamma(\frac{1}{2}+\frac{1}{2}a+\frac{1}{2}b)\Gamma(\frac{1}{2}-\frac{1}{2}a-\frac{1}{2}b+c)}
{\Gamma(\frac{1}{2}+\frac{1}{2}a)\Gamma(\frac{1}{2}+\frac{1}{2}b)\Gamma(\frac{1}{2}-\frac{1}{2}a+c)\Gamma(\frac{1}{2}-\frac{1}{2}b+c)}.\label{2-11}
\emn
\end{thm}
Making the substitutions $a\rightarrow a+1$ and $c\rightarrow c+\frac{1}{2}$ in \eqref{2-11}, we obtain
\bmn
{_3F_2}\ffnk{cc}{1}{1+a,b,\frac{1}{2}+c}{1+\frac{1}{2}a+\frac{1}{2}b,1+2c}
=\frac{\Gamma(\frac{1}{2})\Gamma(1+c)\Gamma(1+\frac{1}{2}a+\frac{1}{2}b)\Gamma(\frac{1}{2}-\frac{1}{2}a-\frac{1}{2}b+c)}
{\Gamma(1+\frac{1}{2}a)\Gamma(\frac{1}{2}+\frac{1}{2}b)\Gamma(\frac{1}{2}-\frac{1}{2}a+c)\Gamma(1-\frac{1}{2}b+c)}.\label{G5}
\emn
By the definition of hypergeometric series and applying \eqref{1.1} and \eqref{1.2}, we can restate \eqref{G5} as
\bmn
\+\+\quad1+b\sum\limits_{k=1}^\infty\frac{\binom{2k}{k}}{k2^{2k}}
\frac{\prod\limits_{i=1}^{k}(1+\frac{a}{i})(1+\frac{2c}{2i-1})\prod\limits_{m=1}^{k-1}(1+\frac{b}{m})}
{\prod\limits_{j=1}^{k}(1+\frac{a+b}{2j})(1+\frac{2c}{j})}.\label{c40}\\
\+\+\quad=\exp\Big\{\sum\limits_{k=1}^\infty(-1)^k\frac{\tau_k}{k}\big[(c-\frac{a}{2}-\frac{b}{2})^k-(\frac{b}{2})^k-(c-\frac{a}{2})^k]\nnm\\
                 \+\+\quad\quad+{\sum\limits_{k=1}^\infty(-1)^k\frac{\sigma_k}{k}\big[c^k+(\frac{a}{2}+\frac{b}{2})^k-(\frac{a}{2})^k-(c-\frac{b}{2})^k]}\Big\}\\
\+\+\quad=\exp\Big\{\tau_1b
+\tau_2(\frac{1}{4}ab-\frac{1}{2}bc)
+\tau_3(\frac{1}{8}a^2b+\frac{1}{8}ab^2-\frac{1}{2}abc-\frac{1}{4}b^2c+\frac{1}{2}bc^2+\frac{1}{12}b^3)
\nnm\\
\+\+\quad\quad+\tau_4(\frac{1}{16}a^3b+\frac{3}{32}a^2b^2-\frac{3}{8}a^2bc+\frac{1}{16}ab^3+\frac{3}{4}abc^2-\frac{3}{8}ab^2c-\frac{1}{8}b^3c+\frac{3}{8}b^2c^2-\frac{1}{2}bc^3)
\nnm\\
\+\+\quad\quad+\tau_5(\frac{1}{32}a^4b+\frac{1}{16}a^3b^2-\frac{1}{4}a^3bc-\frac{3}{8}a^2b^2c+\frac{1}{16}a^2b^3+\frac{3}{4}a^2bc^2+\frac{3}{4}ab^2c^2-abc^3\nnm\\
\+\+\quad\quad-\frac{1}{4}ab^3c+\frac{1}{32}ab^4+\frac{1}{80}b^5-\frac{1}{16}b^4c+\frac{1}{4}b^3c^2-\frac{1}{2}b^2c^3+\frac{1}{2}bc^4)
+\cdots
\Big\}\times\exp\Big\{-\sigma_1b\nnm\\
\+\+\quad\quad
+\sigma_2(\frac{1}{4}ab+\frac{1}{2}bc)+\sigma_3(-\frac{1}{8}a^2b-\frac{1}{8}ab^2-\frac{1}{12}b^3+\frac{1}{4}b^2c-\frac{1}{2}bc^2)
+\sigma_4(\frac{1}{16}a^3b+\frac{3}{32}a^2b^2\nnm\\
\+\+\quad\quad
+\frac{1}{16}ab^3-\frac{3}{8}b^2c^2+\frac{1}{2}bc^3+\frac{1}{8}b^3c)
+\sigma_5(-\frac{1}{32}a^4b-\frac{1}{16}a^3b^2-\frac{1}{16}a^2b^3
-\frac{1}{32}ab^4
+\frac{1}{16}b^4c\nnm\\
\+\+\quad\quad-\frac{1}{4}b^3c^2+\frac{1}{2}b^2c^3-\frac{1}{2}bc^4-\frac{1}{80}b^5)+\cdots
\Big\}.\label{G22}
\emn
Its power series expansion via \eqref{2}--\eqref{5} leads us to infinite summation formula involving generalized harmonic numbers  related to Riemann-Zeta function.


\begin{prop}[Infinite summation formulas related to $\zeta(3)$]
\bmn
\+\+[a^2b]\quad\sum\limits_{k=1}^\infty\frac{\binom{2k}{k}(H_k^{2}-3H_k^{(2)})}{k2^{2k}}=6\zeta(3);\label{g1.8}\\
\+\+[bc^2]\quad\sum\limits_{k=1}^\infty\frac{\binom{2k}{k}[(H_k-O_k)^2+(H_k^{(2)}-O_k^{(2)})]}{k2^{2k}}=\frac{3}{2}\zeta(3);\\
\+\+[abc]\quad\sum\limits_{k=1}^\infty\frac{\binom{2k}{k}H_k(H_k-O_k)}{k2^{2k}}=\frac{7}{2}\zeta(3);\label{G2.4}\\
\+\+[b^3]\quad\sum\limits_{k=1}^\infty\frac{\binom{2k}{k}[(2H_{k-1}-H_k)^2+(H_k^{(2)}-4H_{k-1}^{(2)})]}{k2^{2k}}=4\zeta(3)+\frac{32}{3}\ln^32.
\emn
\end{prop}
Also, we obtain the following infinite summation formulas.
\bmn
\+\+\eqref{g1.6}-\eqref{g1.1}\quad\sum\limits_{k=1}^\infty\frac{\binom{2k}{k}H_k}{k^22^{2k}}=\frac{9}{2}\zeta(3)-\frac{2\pi^2}{3}\ln2;\label{G2.6}\\
\+\+\eqref{g1.6}+\eqref{G2.6}\quad\sum\limits_{k=1}^\infty\frac{\binom{2k}{k}H_k^2}{k2^{2k}}=\frac{21}{2}\zeta(3);\nnm\\
\+\+\eqref{g1.11}+\eqref{g1.7}\quad\sum\limits_{k=1}^\infty\frac{\binom{2k}{k}H_{k-1}}{k^22^{2k}}=\frac{5}{2}\zeta(3)-\frac{\pi^2}{3}\ln2-\frac{4}{3}\ln^32;\label{G2.5}\\
\+\+\eqref{g1.6}-\eqref{G2.5}\quad\quad\sum\limits_{k=1}^\infty\frac{\binom{2k}{k}H_{k-1}^2}{k2^{2k}}=\frac{7}{2}\zeta(3)+\pi^2(\ln2)+\frac{4}{3}\ln^32.\nnm
\emn
\begin{prop}[Infinite summation formulas related to $\zeta(4)$]
\bmn
\+\+[a^2bc]\quad\sum\limits_{k=1}^\infty\frac{\binom{2k}{k}(H_{k}-O_k)(H_k^2-3H_k^{(2)})}{k2^{2k}}=\frac{\pi^4}{4};\label{G2.8}\\
\+\+[abc^2]\quad\sum\limits_{k=1}^\infty\frac{\binom{2k}{k}H_k[(H_{k}-O_k)^2+(H_k^{(2)}-O_k^{(2)})]}{k2^{2k}}=\frac{\pi^4}{8};\label{G2.10}\\
\+\+[a^3b]\quad\sum\limits_{k=1}^\infty\frac{\binom{2k}{k}(H_{k}^3-9{H_k}H_k^{(2)}+14H_k^{(3)})}{k2^{2k}}=\frac{8\pi^4}{15}\label{G2.12}.
\emn
\end{prop}
We also get the following infinite summation formulas related Riemann-zeta functions:
\bnm
\+\+3\times\eqref{g1.20}+\eqref{G2.12}\quad\sum\limits_{k=1}^\infty\frac{\binom{2k}{k}H_k^3+5H_k^{(3)}}{k2^{2k}}=\frac{5\pi^4}{6};\\
\+\+7\times\eqref{g1.20}-\eqref{G2.12}\quad\sum\limits_{k=1}^\infty\frac{\binom{2k}{k}H_k^3+5{H_k}H_k^{(2)}}{k2^{2k}}=\pi^4;\\
\+\+\eqref{g1.20}-\eqref{G2.12}\quad\quad\quad\sum\limits_{k=1}^\infty\frac{\binom{2k}{k}{H_k}H_k^{(2)}-H_k^{(3)}}{k2^{2k}}=\frac{\pi^4}{30}.
\enm
\begin{prop}[Infinite summation formulas related to $\zeta(2)$ and $\zeta(3)$]
\bmn
\+\+[ab^2]\quad\sum\limits_{k=1}^\infty\frac{\binom{2k}{k}(H_k^{(2)}+2{H_k}H_{k-1}-H_k^2)}{k2^{2k}}=3\zeta(3)+\frac{4\pi^2}{3}(\ln2);\label{G2.11}\\
\+\+[b^2c]\quad\sum\limits_{k=1}^\infty\frac{\binom{2k}{k}(H_{k}-O_k)(2H_{k-1}-H_k)}{k2^{2k}}=\frac{3}{2}\zeta(3)+\frac{\pi^2}{3}(\ln2).\label{G3}
\emn
\end{prop}
From the above formulas, we can obtain
\bnm
\+\+\eqref{g1.1.2}-\eqref{g1.1.1}\quad\sum\limits_{k=1}^\infty\frac{\binom{2k}{k}(H_k-O_k)}{k2^{2k}}=\frac{\pi^2}{12};\\
\+\+\eqref{G3}-\eqref{G2.4}\quad\sum\limits_{k=1}^\infty\frac{\binom{2k}{k}(H_k-O_k)}{k^22^{2k}}=\zeta(3)-\frac{\pi^2}{6}(\ln2).
\enm
\begin{prop}[Infinite summation formulas related to $\zeta(2)$ $\zeta(3)$ and $\zeta(4)$]
\bmn
\+\+[a^2b^2]\quad\sum\limits_{k=1}^\infty\frac{\binom{2k}{k}[{H_k^2}H_{k-1}+{H_k^{(2)}}(4H_k-3H_{k-1})]}{k2^{2k}}=\frac{11\pi^4}{45}+12(\ln2)\zeta(3);\label{G2.11}\\
\+\+[b^2c^2]\quad\sum\limits_{k=1}^\infty\frac{\binom{2k}{k}(2H_{k-1}-H_k)[(H_k-O_k)^2+(H_k^{(2)}-O_k^{(2)})]}{k2^{2k}}
=\frac{13\pi^4}{180}+6(\ln2)\zeta(3);\label{G2.9}\\
\+\+[ab^2c]\quad\sum\limits_{k=1}^\infty\frac{\binom{2k}{k}(H_k-O_k)(H_k^{(2)}+2{H_k}H_{k-1}-H_k^2)}{k2^{2k}}
=\frac{13\pi^4}{72}+14(\ln2)\zeta(3).
\emn
\end{prop}
Also, we arrive at the following infinite summation formulas.
\bnm
\+\+\eqref{G2.10}+\eqref{G2.9}\quad\sum\limits_{k=1}^\infty\frac{\binom{2k}{k}H_{k-1}[(H_k-O_k)^2+(H_k^{(2)}-O_k^{(2)})]}{k2^{2k}}
=\frac{71\pi^4}{720}+3(\ln2)\zeta(3);\\
\+\+\eqref{G2.10}-\eqref{G2.9}\quad\sum\limits_{k=1}^\infty\frac{\binom{2k}{k}[(H_k-O_k)^2+(H_k^{(2)}-O_k^{(2)})]}{k^22^{2k}}
=\frac{19\pi^4}{720}-3(\ln2)\zeta(3).
\enm
\section{Summation Formulas Involving  Riemann-Zeta function from Bailey summation theorem}
%
In this section, we shall establish some infinite summation formulas involving Riemann-Zeta function by Bailey's summation theorem \eqref{1.8} and another summation formula \eqref{G2} with two new patterns as follows:
\bnm
\sum_{k=1}^\infty\frac{P_k}{k^i2^k},  \quad i=1,2,  \quad\quad\quad \sum_{k=1}^\infty\frac{3^k}{k^2\binom{2k}{k}}P_k ,
\enm
where $P_k$ is a polynomial in $H_k^{(r)}$ or $O_k^{(r)}$$\:{(k,r\in Z^+)}$.\\

\begin{thm}[Bailey, \cito{Slater}] For complex parameters $a,c$ with $\mathcal{R}(c-1)>0$, the following summation formula is true.
\bmn
{_2F_1}\ffnk{cc}{\frac{1}{2}}{a,1-a}{c}
=\frac{\Gamma(\frac{1}{2}c)\Gamma(\frac{1}{2}+\frac{1}{2}c)}{\Gamma(\frac{1}{2}c+\frac{1}{2}a)\Gamma(\frac{1}{2}+\frac{1}{2}c-\frac{1}{2}a)}.\label{1.8}
\emn
\end{thm}
Performing the replacement $c\rightarrow c+1$ in \eqref{1.8}, we obtain the following expression:
\bmn
{_2F_1}\ffnk{cc}{\frac{1}{2}}{a,1-a}{1+c}
=\frac{\Gamma(\frac{1}{2}+\frac{1}{2}c)\Gamma(1+\frac{1}{2}c)}{\Gamma(\frac{1}{2}+\frac{1}{2}c+\frac{1}{2}a)\Gamma(1+\frac{1}{2}c-\frac{1}{2}a)}.\label{c40}
\emn
Similar to the process illustrated in Theorem 2.1, and after some simplification, some summations involving generalized harmonic numbers related to Riemann-Zeta function can be derived from this identity.



\begin{prop}[Infinite summation formula related to $\zeta(2)$]
\bmn
\+\+[ac]\quad\sum\limits_{k=1}^\infty\frac{H_k}{k2^k}=\frac{\pi^2}{12}.\label{W1.4}
\emn
\end{prop}
\begin{prop}[Infinite summation formula related to $\zeta(3)$]\label{e5}
\bmn
\+\+[ac^2]\quad\sum\limits_{k=1}^\infty\frac{H_k^2+H_k^{(2)}}{k2^k}=\frac{3}{2}\zeta(3).
\emn
\end{prop}
\begin{prop}[Infinite summation formulas related to $\zeta(2)$ and $\zeta(3)$]
\bmn
\+\+[a^2c]\quad\sum\limits_{k=1}^\infty\frac{H_k}{{k^2}2^k}=\zeta(3)-\frac{\pi^2}{12}(\ln2);\\
\+\+[a^3]\quad\sum\limits_{k=1}^\infty\frac{H_{k-1}^{(2)}}{k2^k}=-\frac{1}{4}\zeta(3)+\frac{\pi^2}{12}(\ln2)-\frac{1}{6}\ln^32.
\emn
\end{prop}
\begin{prop}[Infinite summation formula related to $\zeta(4)$]
\bmn
\+\+[ac^3]\quad\sum\limits_{k=1}^\infty\frac{H_k^3+3{H_k}H_k^{(2)}+2H_k^{(3)}}{k2^k}=\frac{7\pi^4}{120}.
\emn
\end{prop}
\begin{prop}[Infinite summation formulas related to $\zeta(2)$, $\zeta(3)$, $\zeta(4)$]
\bmn
\+\+[a^2c^2]\quad\sum\limits_{k=1}^\infty\frac{H_k^2+H_k^{(2)}}{{k^2}2^k}=\frac{19\pi^4}{720}-\frac{3}{2}(\ln2)\zeta(3);\\
\+\+[a^3c]\quad\sum\limits_{k=1}^\infty\frac{{H_k}H_{k-1}^{(2)}}{k2^k}=-\frac{\pi^4}{360}+(\ln2)\zeta(3)-\frac{\pi^2}{24}(\ln^22);\label{W1.2}\\
\+\+[a^4]\quad\sum\limits_{k=1}^\infty\frac{H_{k-1}^{(2)}}{k^22^k}=\frac{\pi^4}{1440}+\frac{1}{4}(\ln2)\zeta(3)-\frac{\pi^2}{24}(\ln^22)+\frac{1}{24}(\ln^42).\label{W1.3}
\emn
\end{prop}
Therefore, we get the following infinite summation formula.
\bnm
\+\+\eqref{W1.2}-\eqref{W1.3}\quad\sum\limits_{k=1}^\infty\frac{{H_{k-1}}H_{k-1}^{(2)}}{k2^k}=-\frac{\pi^4}{288}+\frac{3}{4}(\ln2)\zeta(3)-\frac{\pi^2}{24}(\ln^22).
\enm
\begin{prop}[Infinite summation formulas related to $\zeta(2)$, $\zeta(3)$, $\zeta(4)$, $\zeta(5)$]\label{e5}
\bmn
\+\+[a^2c^3]\quad\sum\limits_{k=1}^\infty\frac{H_k^3+3{H_k}H_k^{(2)}+2H_k^{(3)}}{k^22^k}=12\zeta(5)-\frac{7\pi^4}{120}(\ln2)-\frac{3\pi^2}{8}\zeta(3);\label{W1.1}\\
\+\+[a^3c^2]\quad\sum\limits_{k=1}^\infty\frac{(H_k^2+H_k^{(2)})H_{k-1}^{(2)}}{k2^k}
=-\frac{15}{4}\zeta(5)+\frac{7\pi^4}{240}(\ln2)+\frac{7\pi^2}{24}\zeta(3)-\frac{3}{4}(\ln^22)\zeta(3);\\
\+\+[a^4c]\quad\sum\limits_{k=1}^\infty\frac{{H_k}H_{k-1}^{(2)}}{k^22^k}=-\zeta(5)+\frac{\pi^4}{360}(\ln2)
+\frac{5\pi^2}{48}\zeta(3)-
\frac{1}{2}(\ln^22)\zeta(3)+\frac{\pi^2}{72}(\ln^32).
\emn
\end{prop}
Next, we get some summation formulas related to Riemann-Zeta function due to another summation theorem, which reads as follows.

\begin{thm}[summation theorem \cito{Wang}] For complex parameters $b,d$ with $\mathcal{R}(2b+2d-3)<0$, the following summation formula is true.
\bmn
{_3F_2}\ffnk{cc}{\frac{3}{4}}{b,d,\frac{b+d}{3}}{\frac{b+d}{2},\frac{1+b+d}{2}}
=\frac{\Gamma(1+b+d)\Gamma(1+\frac{1}{3}b)\Gamma(1+\frac{1}{3}d)}{\Gamma(1+b)\Gamma(1+d)\Gamma(1+\frac{b+d}{3})}.\label{G2}
\emn
\end{thm}

By the definition of hypergeometric series and the identities \eqref{1.1} and \eqref{1.2}, the summation formula \eqref{G2} can be reformulated as
\bmn
\+\+\quad1+\frac{2}{3}{bd}\sum\limits_{k=1}^\infty\frac{3^k\prod\limits_{i=1}^{k-1}(1+\frac{b}{i})
(1+\frac{d}{i})(1+\frac{b+d}{3i})}{k^2\binom{2k}{k}\prod\limits_{j=1}^{k-1}(1+\frac{b+d}{2j})\prod\limits_{m=1}^{k}(1+\frac{b+d}{2m-1})}\\
\+\+\quad=\exp\Big\{\sum\limits_{k=1}^\infty\frac{(-1)^k\sigma_k}{k}\big[(b+d)^k+(\frac{b}{3})^k+(\frac{d}{3})^k-b^k-d^k-(\frac{b+d}{3})^k\big]\Big\}\nnm\\
\+\+\quad=\exp\Big\{\sigma_2(\frac{8}{9}bd)+\sigma_3(-\frac{26}{27}b^2d-\frac{26}{27}bd^2)
+\sigma_4(\frac{80}{81}b^3d+\frac{80}{81}bd^3+\frac{120}{81}b^2d^2)\nnm\\
\+\+\quad\quad+\sigma_5(-\frac{242}{243}b^4d-\frac{484}{243}b^3d^2-\frac{484}{243}b^2d^3-\frac{242}{243}bd^4)+\cdots
\Big\}.\label{G22}
\emn
Some infinite summation formulas involving generalized harmonic numbers related to Riemann-Zeta function obtained from its power series expansion are as follows.

\begin{prop}[Infinite summation formula related to $\zeta(3)$]
\bmn
\+\+[b^2d]\quad\quad\:\:\sum\limits_{k=1}^\infty\frac{3^k(6O_k-5H_{k-1})}{k^2\binom{2k}{k}}=\frac{182}{3}\zeta(3).\label{W1.1}
\emn
\end{prop}

\begin{prop}[Infinite summation formula related to $\zeta(2)$, $\zeta(4)$]
\bmn
\+\+[b^2d^2]\quad\quad\:\:\sum\limits_{k=1}^\infty\frac{3^k[36(O_k^2+O_k^{(2)})+25H_{k-1}^2+5H_{k-1}^{(2)}-60{H_{k-1}}O_k]}{k^2\binom{2k}{k}}=\frac{56\pi^4}{3}.
\emn
\end{prop}

%


There are many other infinite summation formulas involving generalized harmonic numbers related to Riemann-Zeta function can be established from these summation formulas. Here, we just present some results for examples, the interested reader can do by themselves.

\end{document}